# On the convergence of periodic Navier-Stokes flows

David T. Purvance

**Abstract**

The 3D spatially periodic Navier-Stokes equation is posed as a nonlinear matrix differential equation. When the flow is assumed to be a time series having unknown wavenumber coefficients, then the matrix in this periodic Navier-Stokes matrix differential equation becomes a time series of matrices. Posed in this way, any flow's unknown wavenumber coefficients can be solved for recursively beginning with the zeroth-order coefficient representing any initial flow. When all matrices in the time series commute, a flow's unknown coefficients also represent the Taylor expansion of a stable matrix exponential product operating on the initial flow. This paper argues that a solution's coefficients converge to these bounded Taylor coefficients when these bounds are evaluated with a general solution's noncommutative matrices.

## 1.0 Introducing the Periodic Navier-Stokes Matrix Differential Equation

Begin with the $\mathbb{C}^{3 \times 1}$ Navier-Stokes vector equation [1] describing the time evolution of an incompressible flow's $M$ Fourier modes $\mathbf{u}(\boldsymbol{\kappa}_j, t)$

$$\left(\frac{d}{dt} + \nu \kappa_j^2\right) \mathbf{u}(\boldsymbol{\kappa}_j, t) = -P(\boldsymbol{\kappa}_j) \sum_{k=1}^{M} J(\boldsymbol{\kappa}_j, \boldsymbol{\kappa}'_k, t) \mathbf{u}(\boldsymbol{\kappa}'_k, t) \qquad (1)$$

where $\boldsymbol{\kappa}_j$ and $\boldsymbol{\kappa}'_k$ are $M = (L+1)^3$ discrete wavenumbers $(l_1 d\kappa, l_2 d\kappa, l_3 d\kappa)$ with $l_1, l_2, l_3 = (-L/2, ..., 0, ..., L/2)$ for some even integer $L$ and where $d\kappa = 2\pi/L$. Included in (1) are $\kappa_j^2 = |\boldsymbol{\kappa}_j|^2$, a known viscosity $\nu > 0$, $P(\boldsymbol{\kappa}_j)$ the projection tensor

$$P(\boldsymbol{\kappa}_j) = \begin{bmatrix} 1 - \dfrac{\kappa_{j,1}^2}{\kappa_j^2} & -\dfrac{\kappa_{j,1}\kappa_{j,2}}{\kappa_j^2} & -\dfrac{\kappa_{j,1}\kappa_{j,3}}{\kappa_j^2} \\ -\dfrac{\kappa_{j,2}\kappa_{j,1}}{\kappa_j^2} & 1 - \dfrac{\kappa_{j,2}^2}{\kappa_j^2} & -\dfrac{\kappa_{j,2}\kappa_{j,3}}{\kappa_j^2} \\ -\dfrac{\kappa_{j,3}\kappa_{j,1}}{\kappa_j^2} & -\dfrac{\kappa_{j,3}\kappa_{j,2}}{\kappa_j^2} & 1 - \dfrac{\kappa_{j,3}^2}{\kappa_j^2} \end{bmatrix}, \qquad (2)$$

$$J(\boldsymbol{\kappa}_j, \boldsymbol{\kappa}'_k, t) = \begin{bmatrix} i\kappa_{j,1}u_1(\boldsymbol{\kappa}_j - \boldsymbol{\kappa}'_k, t) & i\kappa_{j,2}u_1(\boldsymbol{\kappa}_j - \boldsymbol{\kappa}'_k, t) & i\kappa_{j,3}u_1(\boldsymbol{\kappa}_j - \boldsymbol{\kappa}'_k, t) \\ i\kappa_{j,1}u_2(\boldsymbol{\kappa}_j - \boldsymbol{\kappa}'_k, t) & i\kappa_{j,2}u_2(\boldsymbol{\kappa}_j - \boldsymbol{\kappa}'_k, t) & i\kappa_{j,3}u_2(\boldsymbol{\kappa}_j - \boldsymbol{\kappa}'_k, t) \\ i\kappa_{j,1}u_3(\boldsymbol{\kappa}_j - \boldsymbol{\kappa}'_k, t) & i\kappa_{j,2}u_3(\boldsymbol{\kappa}_j - \boldsymbol{\kappa}'_k, t) & i\kappa_{j,3}u_3(\boldsymbol{\kappa}_j - \boldsymbol{\kappa}'_k, t) \end{bmatrix}, \qquad (3)$$

and the implied incompressibility condition $\boldsymbol{\kappa}_j \cdot \mathbf{u}(\boldsymbol{\kappa}_j, t) = 0$. Assume a known, finite initial flow $\mathbf{u}(\boldsymbol{\kappa}_j, t=0) = \mathbf{u}_0(\boldsymbol{\kappa}_j)$.

Simultaneously consider all wavenumbers $\boldsymbol{\kappa}_j$ by writing out (1) in consecutive trios of rows for $j = 1, \ldots, M$. This forms the $\mathbb{C}^{3M \times 3M}$ periodic Navier-Stokes matrix differential equation

$$\frac{d\mathbf{u}}{dt} = U(\mathbf{u})\mathbf{u} \qquad (4)$$

where

$$\mathbf{u} = \begin{bmatrix} \mathbf{u}(\boldsymbol{\kappa}_1, t) \\ \vdots \\ \mathbf{u}(\boldsymbol{\kappa}_M, t) \end{bmatrix} = \begin{bmatrix} \mathbf{u}(\boldsymbol{\kappa}'_1, t) \\ \vdots \\ \mathbf{u}(\boldsymbol{\kappa}'_M, t) \end{bmatrix} \qquad (5)$$

and

$$U(\mathbf{u}) = D + PJ(\mathbf{u}), \qquad (6)$$

with

$$D = \begin{bmatrix} -D(\kappa_1^2) & \cdots & 0 \\ \vdots & \ddots & \vdots \\ 0 & \cdots & -D(\kappa_M^2) \end{bmatrix}, \tag{7}$$

$$D(\kappa_j^2) = \begin{bmatrix} \nu\kappa_j^2 & 0 & 0 \\ 0 & \nu\kappa_j^2 & 0 \\ 0 & 0 & \nu\kappa_j^2 \end{bmatrix}, \tag{8}$$

$$P = \begin{bmatrix} -P(\boldsymbol{\kappa}_1) & \cdots & 0 \\ \vdots & \ddots & \vdots \\ 0 & \cdots & -P(\boldsymbol{\kappa}_M) \end{bmatrix}, \tag{9}$$

and

$$J(\mathbf{u}) = \begin{bmatrix} J(\boldsymbol{\kappa}_1, \boldsymbol{\kappa}_1', t) & \cdots & J(\boldsymbol{\kappa}_1, \boldsymbol{\kappa}_M', t) \\ \vdots & \ddots & \vdots \\ J(\boldsymbol{\kappa}_M, \boldsymbol{\kappa}_1', t) & \cdots & J(\boldsymbol{\kappa}_M, \boldsymbol{\kappa}_M', t) \end{bmatrix}. \tag{10}$$

Equation (4) becomes the continuous periodic Navier-Stokes equations as $M \to \infty$.

## 2.0 Series Solution to the Navier-Stokes Matrix Differential Equation

First, expand the wavenumber-shifted portion of $\mathbf{u}(\boldsymbol{\kappa}_j, t)$ in $J(\boldsymbol{\kappa}_j, \boldsymbol{\kappa}_k', t)$, i.e., in $J(\boldsymbol{\kappa}_j, \boldsymbol{\kappa}_k', t)$ let

$$\mathbf{u}(\boldsymbol{\kappa}_j - \boldsymbol{\kappa}_k', t) = \sum_{n=0}^{\infty} \mathbf{u}_n(\boldsymbol{\kappa}_j - \boldsymbol{\kappa}_k') t^n \tag{11}$$

for unknown coefficients $\mathbf{u}_n(\boldsymbol{\kappa}_j)$ apart from the initial $\mathbf{u}_0(\boldsymbol{\kappa}_j)$. This creates the equivalent expression of (1)

$$\left(\frac{d}{dt}+\nu\kappa_j^2\right)\mathbf{u}(\mathbf{\kappa}_j,t)=-\sum_{n=0}^{\infty}\left(P(\mathbf{\kappa}_j)\sum_{k=1}^{M}J_n(\mathbf{\kappa}_j,\mathbf{\kappa}_k')\mathbf{u}(\mathbf{\kappa}_k',t)\right)t^n \qquad (12)$$

with

$$J_n(\mathbf{\kappa}_j,\mathbf{\kappa}_k')=\begin{bmatrix} i\kappa_{j,1}u_{n,1}(\mathbf{\kappa}_j-\mathbf{\kappa}_k') & i\kappa_{j,2}u_{n,1}(\mathbf{\kappa}_j-\mathbf{\kappa}_k') & i\kappa_{j,3}u_{n,1}(\mathbf{\kappa}_j-\mathbf{\kappa}_k') \\ i\kappa_{j,1}u_{n,2}(\mathbf{\kappa}_j-\mathbf{\kappa}_k') & i\kappa_{j,2}u_{n,2}(\mathbf{\kappa}_j-\mathbf{\kappa}_k') & i\kappa_{j,3}u_{n,2}(\mathbf{\kappa}_j-\mathbf{\kappa}_k') \\ i\kappa_{j,1}u_{n,3}(\mathbf{\kappa}_j-\mathbf{\kappa}_k') & i\kappa_{j,2}u_{n,3}(\mathbf{\kappa}_j-\mathbf{\kappa}_k') & i\kappa_{j,3}u_{n,3}(\mathbf{\kappa}_j-\mathbf{\kappa}_k') \end{bmatrix}. \qquad (13)$$

In $3M$ space (11), (12) and (13) rewrite (4) as

$$\frac{d\mathbf{u}}{dt}=\left(\sum_{n=0}^{\infty}U_n(\mathbf{u}_n)t^n\right)\mathbf{u} \qquad (14)$$

with

$$U_n=\delta_{0,n}D+PJ_n. \qquad (15)$$

In contrast to $U(\mathbf{u})$, $U_n(\mathbf{u}_n)$ are independent of time.

Finally, use unshifted (11) to rewrite (14) as

$$\sum_{n=0}^{\infty}n\mathbf{u}_n t^{n-1}=\left(\sum_{p=0}^{\infty}U_p t^p\right)\sum_{q=0}^{\infty}\mathbf{u}_q t^q=\sum_{p=0}^{\infty}\sum_{q=0}^{\infty}U_p\mathbf{u}_q t^{p+q}. \qquad (16)$$

Matching coefficients in (16) when $t^{n-1}=t^{p+q}$ recursively solves for the unknown flow coefficients $\mathbf{u}_n$. For $n>0$ they are

$$\mathbf{u}_n=\frac{1}{n}\sum_{p=0}^{n-1}U_p\mathbf{u}_{n-1-p}. \qquad (17)$$

When expanded and simplified, the first few orders of (17), generated by Maple® 11's Physics Package using noncommutative $U_n$, are

$$\mathbf{u}_1 = U_0 \mathbf{u}_0 = S_1 \mathbf{u}_0$$

$$\mathbf{u}_2 = \frac{1}{2}\left(U_0^2 + U_1\right)\mathbf{u}_0 = S_2 \mathbf{u}_0$$

$$\mathbf{u}_3 = \frac{1}{24}\left(4U_0^3 + 4U_0 U_1 + 8U_1 U_0 + 8U_2\right)\mathbf{u}_0 = S_3 \mathbf{u}_0 \quad (18)$$

$$\mathbf{u}_4 = \frac{1}{24}\left(U_0^4 + U_0^2 U_1 + 2U_0 U_1 U_0 + 2U_0 U_2 + 3U_1 U_0^2 + 3U_1^2 + 6U_2 U_0 + 6U_3\right)\mathbf{u}_0 = S_4 \mathbf{u}_0$$

$$\vdots$$

Note that $\mathbf{u}_n$ consist of various matrix product terms which when summed, transform the initial flow coefficient $\mathbf{u}_0$ into the $n^{\text{th}}$ order solution coefficient, i.e.,

$$\mathbf{u}_n = S_n \mathbf{u}_0, \quad (19)$$

and

$$S_n = \sum_p s_{np} \prod_{g \in G_{np}} U_g \quad (20)$$

with rational fractions $s_{np}$ and $G_{np}$ the set of time-order indices describing the left-to-right matrix product order of $\prod U_g$. For instance, from (18) $s_{43} = 1/12$ and $G_{43}$ is the set of time-order indices $\{0,1,0\}$. Simple bookkeeping argues that the number of terms making up $S_n$ and counted by index $p(n)$ grows by $2^{n-1}$.

### 3.0 A Convergent Function of Noncommutative $U_n$

When $U_n$ commute, the solution to the Navier-Stokes (14) is

$$\mathbf{u} = \exp\left(\sum_{n=0}^{\infty} \frac{1}{n+1} U_n t^{n+1}\right)\mathbf{u}_0 = \prod_{n=0}^{\infty} \exp\left(\frac{1}{n+1} U_n t^{n+1}\right)\mathbf{u}_0. \quad (21)$$

The, matrix-exponential-product (far right) version of this solution is convergent irrespective of the commutative properties of $U_n$ because individually the Navier-Stokes $U_n$ are inherently stable, proven in the Appendix.

Let function **b** be the Taylor expansion the matrix exponential products in (21) using a general solution's noncommutative $U_n$

$$\mathbf{b} = \sum_{n=0}^{\infty} \mathbf{b}_n t^n = \prod_{n=0}^{\infty} \left( \sum_{p=0}^{\infty} \frac{1}{p!} \left( \frac{1}{n+1} U_n t^{n+1} \right)^p \right) \mathbf{u}_0 . \tag{22}$$

While time series **b** converges, it obviously is not a solution to the Navier-Stokes (14) because, rather than commutative $U_n$, a general solution's noncommutative $U_n$ are assumed.

Again using Maple® 11, expansion (22) using noncommutative $U_n$ produces coefficients $\mathbf{b}_n$ whose first few terms are

$$\begin{aligned}
\mathbf{b}_1 &= U_0 \mathbf{u}_0 = B_1 \mathbf{u}_0 \\
\mathbf{b}_2 &= \frac{1}{2} \left( U_0^2 + U_1 \right) \mathbf{u}_0 = B_2 \mathbf{u}_0 \\
\mathbf{b}_3 &= \frac{1}{24} \left( 4 U_0^3 + 12 U_0 U_1 + 8 U_2 \right) \mathbf{u}_0 = B_3 \mathbf{u}_0 \\
\mathbf{b}_4 &= \frac{1}{24} \left( U_0^4 + 6 U_0^2 U_1 + 8 U_0 U_2 + 3 U_1^2 + 6 U_3 \right) \mathbf{u}_0 = B_4 \mathbf{u}_0 \\
&\vdots
\end{aligned} \tag{23}$$

Like $\mathbf{u}_n$, $\mathbf{b}_n$ consist of various matrix products terms, made from noncommutative $U_0,...,U_{n-1}$, which when summed, transform the initial flow coefficient $\mathbf{u}_0$ into the $n^{th}$ order bound $\mathbf{b}_n$ with

$$\mathbf{b}_n = B_n \mathbf{u}_0 \qquad (24)$$

where

$$B_n = \sum_q b_{nq} \prod_{g \in G_{nq}} U_g \qquad (25)$$

with rational fractions $b_{nq}$ and $G_{nq}$ the matrix product order of $\prod U_g$. Because **b** is stable, $B_n$ are bounded and converge to zero.

## 4.0 Convergence

Form the difference between $S_n$ which is not known to converge and $B_n$ which is known to converge

$$S_n - B_n = \sum_p s_{np} \prod_{g \in G_{np}} U_g - \sum_q b_{nq} \prod_{g \in G_{nq}} U_g = \sum_r d_{nr} \prod_{g \in G_{nr}} U_g. \qquad (26)$$

When matrix product order is immaterial, corresponding to commutative $U_n$, then

$$S_n - B_n = 0 \Rightarrow \sum_r d_{nr} = 0. \qquad (27)$$

However, since $d_{nr}$ are fixed, (27) is true regardless of the commutativity properties of $U_n$. For example, the first 4 orders of $S_n - B_n$ are

$$\begin{aligned}
\mathbf{u}_1 - \mathbf{b}_1 &= 0 \\
\mathbf{u}_2 - \mathbf{b}_2 &= 0 \\
\mathbf{u}_3 - \mathbf{b}_3 &= \frac{1}{3}(U_1 U_0 - U_0 U_1)\mathbf{u}_0 \\
\mathbf{u}_4 - \mathbf{b}_4 &= \frac{1}{24}(2U_0 U_1 U_0 + 3U_1 U_0^2 - 5U_0^2 U_1)\mathbf{u}_0 + \frac{1}{4}(U_2 U_0 - U_0 U_2)\mathbf{u}_0
\end{aligned} \qquad (28)$$

While one does not see in (28) the nested commutators of quantum mechanics, one does see that differences $S_n - B_n$ do indeed vanish when commutators $[U_0, U_1]$ and $[U_0, U_2]$ vanish, and, that the sum of the fractions in $S_n - B_n$ vanish regardless of the commutativity properties of $U_n$.

Convergence in the limit

$$\lim_{n \to \infty}(S_n - B_n) = \lim_{r(n) \to \infty} \sum_r d_{nr} \prod_{g \in G_{nr}} U_g = \lim_{r(n) \to \infty} \left( \sum_r d_{nr} \right) \left( \sum_r \prod_{g \in G_{nr}} U_g \right) = 0 \qquad (29)$$

is reached if over $r(n)$ $d_{nr}$ and elements of matrices $\prod_{g \in G_{nr}} U_g$ approach linear independence. Since coefficients $d_{nr}$ are the same for all flows, they must be independent of the specifics of any flow $\prod_{g \in G_{nr}} U_g$ which argues that (29) is true.

**5.0 Summary and Conclusion**

In summary the 3D spatially periodic Navier-Stokes equation has been posed as a nonlinear matrix differential equation $d\mathbf{u}/dt = U(\mathbf{u})\mathbf{u}$ of dimension $3M$. When the flow $\mathbf{u}$ is assumed to be a time series of order $N$ with unknown coefficients $\mathbf{u}_n$, then matrix $U$ becomes a time series of matrices $U_n$ of order $N$. $U_n$ are inherently stable. So posed, a flow's unknown coefficients $\mathbf{u}_n = S_n \mathbf{u}_0$ for transformations $S_n$ can be solved for recursively starting with any known initial flow $\mathbf{u}_0$.

When $U_n$ commute, the time series solution remains bounded for all time because coefficients $\mathbf{u}_n$ also represent the Taylor expansion coefficients of the product of stable matrix exponentials $\mathbf{b}_n = B_n \mathbf{u}_0$. Solution transformations $S_n$ and bound transformations $B_n$ come in the form of matrix polynomials $\sum_p s_{np} \prod_{g \in G_{np}} U_g$ and $\sum_q b_{nq} \prod_{g \in G_{nq}} U_g$, respectively, having fixed positive rational fractions $s_{np}$ and $b_{nq}$. In the differences $S_n - B_n = \sum_{r(n)} d_{nr} \prod_{g \in G_{nr}} U_g$ the fractions always sum to zero regardless of the commutativty properties of $U_0, ..., U_{n-1}$, i.e., $\sum_p s_{np} - \sum_q b_{nq} = \sum_r d_{nr} = 0$. Scalars $d_{nr}$ are fixed for all flows and are not a function of the specifics of any flow $\prod_{g \in G_{nr}} U_g$. This paper therefore argues that in the limit these two components making up the inner product $\sum_{r(n)} d_{nr} \prod_{g \in G_{nr}} U_g$ become linearly independent. Consequently, as $n$ gets large, $S_n - B_n \to 0$, and solution transformations $S_n$ converge onto bounded transformations $B_n$.

The conclusion is that, as $M, N \to \infty$, $\mathbf{u}_n$ become the convergent time series solution coefficients to the continuous Navier-Stokes equations.

# Appendix

*$U_n$ Stability*

Recognize that for the inverse Fourier transform of the three spatial components of $\mathbf{u}_n(\boldsymbol{\kappa}_j)$ to be real-valued, then they must be conjugate symmetric, i.e., $\mathbf{u}_n(\boldsymbol{\kappa}_j) = \mathbf{u}_n^*(-\boldsymbol{\kappa}_j)$ where the asterisk superscript denotes complex conjugate. Therefore, for $n > 0$ the Fourier conjugate symmetry of $U_n(\boldsymbol{\kappa}_j, \boldsymbol{\kappa}_k') = P(\boldsymbol{\kappa}_j) J_n(\boldsymbol{\kappa}_j, \boldsymbol{\kappa}_k')$ is

$$\begin{aligned} U_n^*(-\boldsymbol{\kappa}_j, -\boldsymbol{\kappa}_k') &= \left(P(-\boldsymbol{\kappa}_j) J_n(-\boldsymbol{\kappa}_j, -\boldsymbol{\kappa}_k')\right)^* \\ &= P^*(-\boldsymbol{\kappa}_j) J_n^*(-\boldsymbol{\kappa}_j, -\boldsymbol{\kappa}_k') \\ &= P(\boldsymbol{\kappa}_j) J_n(\boldsymbol{\kappa}_j, \boldsymbol{\kappa}_k') \\ &= U_n(\boldsymbol{\kappa}_j, \boldsymbol{\kappa}_k') \end{aligned} \qquad (30)$$

which is true because from (2) $P(\boldsymbol{\kappa}_j)$ is real and symmetric with $P(\boldsymbol{\kappa}_j) = P(-\boldsymbol{\kappa}_j)$ and because from (13) each $(j,l)$ element of $J_n(\boldsymbol{\kappa}_j, \boldsymbol{\kappa}_k')$ has the property

$$(-i\kappa_{j,l})^* u_{n,l}^*(-(\boldsymbol{\kappa}_j - \boldsymbol{\kappa}_k')) = i\kappa_{j,l} u_{n,l}(\boldsymbol{\kappa}_j - \boldsymbol{\kappa}_k'). \qquad (31)$$

This same result is also holds for $U_0$ because $D(\boldsymbol{\kappa}_j) = D(-\boldsymbol{\kappa}_j)$.

It is important to distinguish this Fourier conjugate symmetry from the Hermitian or diagonal conjugate symmetry of $\mathbb{C}^{3\times 3}$ $U_n(\boldsymbol{\kappa}_j, \boldsymbol{\kappa}_k')$ which would be denoted by $U_n(\boldsymbol{\kappa}_j, \boldsymbol{\kappa}_k') = U_n^*(\boldsymbol{\kappa}_k', \boldsymbol{\kappa}_j) = U_n^H(\boldsymbol{\kappa}_j, \boldsymbol{\kappa}_k')$ where the superscript $H$ denotes the Hermitian or matrix conjugate transpose. There is no diagonal conjugate symmetry in $U_n(\boldsymbol{\kappa}_j, \boldsymbol{\kappa}_k')$ because there is no diagonal conjugate symmetry in $J_n(\kappa_j, \kappa'_k)$. Furthermore, a simple inspection the $U_n(\kappa_j, \kappa'_k)$ $j = k$ diagonals of $U_n$ proves that $U_n$ is neither Hermitian nor skew-Hermitian, i.e.,

$$U_n \neq \pm U_n^H. \tag{32}$$

Nonetheless, for $n > 0$ and in terms of submatrix products the eigenvalues of $U_n$ are

$$\lambda_{m,n} = \sum_{j=1}^{M} \mathbf{x}_{m,n}^H(\kappa_j) \sum_{k=1}^{M} U_n(\kappa_j, \kappa'_k) \mathbf{x}_{m,n}(\kappa'_k) \\ = \sum_{j=1}^{M} \mathbf{x}_{m,n}^H(\kappa_j) \sum_{k=1}^{M} U_n(\kappa_j, \kappa'_k) \mathbf{x}_{m,n}(\kappa'_k) \tag{33}$$

where $\mathbf{x}_{m,n}$ is the $m^{\text{th}}$ eigenvector of $U_n$. The complex conjugate of (33) is

$$\lambda_{m,n}^* = \left( \sum_{j=1}^{M} \sum_{k=1}^{M} \mathbf{x}_{m,n}^H(\kappa_j) U_n(\kappa_j, \kappa'_k) \mathbf{x}_{m,n}(\kappa'_k) \right)^* \\ = \sum_{j=1}^{M} \sum_{k=1}^{M} \mathbf{x}_{m,n}^T(\kappa_j) U_n^*(\kappa_j, \kappa'_k) \mathbf{x}_{m,n}^*(\kappa'_k) \tag{34}$$

Equation (34) remains unchanged if in (34) $\kappa_j$ and $\kappa'_k$ are replaced with $-\kappa_j$ and $-\kappa'_k$. This amounts to performing the wavenumber summations in the negative order of (34). Doing so

$$\lambda_{m,n}^* = \sum_{j=1}^{M} \mathbf{x}_{m,n}^T(-\kappa_j) \sum_{k=1}^{M} U_n^*(-\kappa_j, -\kappa'_k) \mathbf{x}_{m,n}^*(-\kappa'_k). \tag{35}$$

However, by (30) the conjugate symmetry of $U_n(\kappa_j, \kappa'_k)$ allows (35) to be rewritten

$$\lambda_{m,n}^* = \sum_{j=1}^{M} \mathbf{x}_{m,n}^T(-\kappa_j) \sum_{k=1}^{M} U_n(\kappa_j, \kappa'_k) \mathbf{x}_{m,n}^*(-\kappa'_k) \tag{36}$$

Comparing (36) to (33) identifies a second eigenvector of $U_n$, denoted by the subscript $l$, that is Fourier conjugate symmetric to $\mathbf{x}_{m,n}$

$$\mathbf{x}_{l,n}(\kappa_j) = \mathbf{x}^*_{m,n}(-\kappa_j) \qquad (37)$$

and a second eigenvalue $\lambda_{l,n}$ conjugate to $\lambda_{m,n}$

$$\lambda^*_{m,n} = \sum_{j=1}^{M} \mathbf{x}^H_{l,n}(\kappa_j) \sum_{k=1}^{M} U_n(\kappa_j, \kappa'_k) \mathbf{x}_{l,n}(\kappa'_k) = \lambda_{l,n}. \qquad (38)$$

Form the sum

$$2(\lambda_{m,n} + \lambda^*_{m,n}) = 2(\lambda_{l,n} + \lambda^*_{l,n}) = \mathbf{x}^H_{m,n}(U_n + U^H_n)\mathbf{x}_{m,n} = \mathbf{x}^H_{l,n}(U^H_n + U_n)\mathbf{x}_{l,n} \qquad (39)$$

Recognizing (32) and in particular that $U_n$ is not skew Hermitian, let

$$H_n = \frac{1}{2}(U_n + U^H_n) \neq 0. \qquad (40)$$

Write (39) as

$$2(\lambda_{m,n} + \lambda^*_{m,n}) = 2(\lambda_{m,n} + \lambda_{l,n}) = 2(\lambda_{l,n} + \lambda^*_{l,n}) = \mathbf{x}^H_{m,n} H_n \mathbf{x}_{m,n} = \mathbf{x}^H_{l,n} H_n \mathbf{x}_{l,n} = 2\mu_{m,n}. \qquad (41)$$

for some yet-to-be-determined $\mu_{m,n}$. From (41)

$$H_n \mathbf{x}_{m,n} = 2\mu_{m,n} \mathbf{x}_{m,n} = H_n \mathbf{x}_{l,n} = 2\mu_{m,n} \mathbf{x}_{l,n} \qquad (42)$$

giving

$$H_n \mathbf{x}_{m,n} - H_n \mathbf{x}_{l,n} = 2\mu_{m,n} \mathbf{x}_{m,n} - 2\mu_{m,n} \mathbf{x}_{l,n} = 2\mu_{m,n}(\mathbf{x}_{m,n} - \mathbf{x}_{l,n}) = 0. \qquad (43)$$

From (37) $(\mathbf{x}_{m,n}, \mathbf{x}_{l,n})$ are Fourier conjugate pairs and therefore

$$\mathbf{x}_{m,n} - \mathbf{x}_{l,n} \neq 0 \qquad (44)$$

which implies in (43)

$$\mu_{m,n} = \lambda_{m,n} + \lambda^*_{m,n} = \lambda_{m,n} + \lambda_{l,n} = \lambda_{l,n} + \lambda^*_{l,n} = 0. \qquad (45)$$

Equation (45) can be true if only if $(\lambda_{m,n}, \lambda^*_{m,n}) = (\lambda_{m,n}, \lambda_{l,n})$ are purely imaginary conjugate pairs, i.e., if and only if

$$\left(\lambda_{m,n}, \lambda_{l,n}\right) = \left(i\omega_{m,n}, -i\omega_{m,n}\right) \tag{46}$$

for real pairs $\pm\omega_{m,n}$.

Since the eigenvalues of diagonal $D$ are real and non-positive ($-\nu\kappa_j^2$) and since the eigenvalues of $PJ_0$ are also zero or purely imaginary by the same (30)-(46) argument, then there are no real positive eigenvalues in the whole $U_n$ lot, rendering them all stable.

**References**

[1] Pope SB 2003 *Turbulent Flows* eq 6.146 Cambridge University Press NY NY

13